\documentclass[a4paper,11pt,reqno]{amsart}
\usepackage[latin1]{inputenc}
\usepackage[english]{babel}
\usepackage{amssymb,amsfonts, amsmath, color}
\usepackage[margin=2.7cm]{geometry}
\usepackage{xcolor}
\usepackage{graphicx, color, enumerate}
\usepackage[latin1]{inputenc}
\usepackage[active]{srcltx}
\usepackage{tikz}
\usepackage{pgf}
\usepackage{etex}
\usepackage{verbatim}
\usepackage{tikz-3dplot}
\usepackage{amsmath}
\usepackage{amsfonts}
\usepackage{amssymb}
\usepackage{amsthm}
\usepackage{algpseudocode}
\usepackage{float}
\usepackage{xcolor}
\usepackage{mathrsfs}
\usepackage{import}
\usepackage{geometry}
\usepackage{fancyhdr}
\usepackage{fp}

\newtheorem{theorem}{Theorem}[section]
\newtheorem{proposition}[theorem]{Proposition}

\newtheorem{lemma}[theorem]{Lemma}
\newtheorem{remark}[theorem]{Remark}
\newtheorem{definition}[theorem]{Definition}

\newcommand{\bcl}{\begin{center}}
\newcommand{\ecl}{\end{center}}
\newcommand{\brl}{\begin{right}}
\newcommand{\erl}{\end{right}}
\newcommand{\ben}{\begin{enumerate}}
\newcommand{\een}{\end{enumerate}}
\newcommand{\overliner}{\begin{array}}
\newcommand{\earr}{\end{array}}
\newcommand{\btab}{\begin{tabular}}
\newcommand{\etab}{\end{tabular}}
\newcommand{\bdoc}{\begin{document}}
\newcommand{\edoc}{\end{document}}
\newcommand{\beqy}{\begin{eqnarray}}
\newcommand{\eeqy}{\end{eqnarray}}

\newcommand{\beqi}{\begin{eqnarray*}}
\newcommand{\eeqi}{\end{eqnarray*}}
\newcommand{\bitem}{\begin{itemize}}

\newcommand{\eitem}{\end{itemize}}
\newcommand{\nln}{\newline}
\newcommand{\newt}{\newtheorem}

\newcommand{\pa}{\partial}
\newcommand{\re}{{I\!\!R}}
\newcommand{\Rn}{\R^N}
\newcommand{\xr}{x\in\R }
\newcommand{\x}{\times}
\newcommand{\dyle}{\displaystyle}
\newcommand{\ene}{{I\!\!N}}
\newcommand{\irn}{\int\limits_{\R^N}}
\newcommand{\io}{\int\limits_{\O}}
\newcommand{\meas}{{\rm meas\,}}
\newcommand{\dif}{\nabla_{xy}}
\newcommand{\sign}{{\rm sign}}
\newcommand{\map}{\longrightarrow }
\newcommand{\imp}{\Longrightarrow }
\renewcommand{\div}{\nabla\cdot }
\newcommand{\sen}{{\rm sen\,}}
\newcommand{\tg}{{\rm tg\,}}
\newcommand{\arcsen}{{\rm arcsen\,}}
\newcommand{\arctg}{{\rm arctg\,}}
\newcommand{\supp}{{\textsl supp\ }}
\newcommand{\ity}{\int_{-\iy}^{+\iy}}
\newcommand{\limit}{\lim\limits}
\newcommand{\limi}{\limit_{n\to\infty}}
\newcommand{\sumi}{\sum\limits_{n=1}^{\infty}}
\newcommand{\ulu}{\underline u}
\newcommand{\ulw}{\underline w}
\newcommand{\ulz}{\underline z}
\newcommand{\ulv}{\underline v}
\newcommand{\uls}{\underline s}
\newcommand{\olu}{\overline u}
\newcommand{\olv}{\overline v}
\newcommand{\ols}{\overline s}
\newcommand{\ob}{\overline\b}
\newcommand{\ovar}{\overline\var}
\newcommand{\wv}{\widetilde v}
\newcommand{\wu}{\widetilde u}
\newcommand{\ws}{\widetilde s}
\renewcommand{\a }{\alpha }
\renewcommand{\b }{\beta }
\newcommand{\g }{\gamma}
\newcommand{\G }{\Gamma }
\renewcommand{\d }{\delta }

\newcommand{\D }{\Delta }
\newcommand{\e }{\varepsilon }
\newcommand{\z }{\zeta }
\renewcommand{\l }{\lambda }
\renewcommand{\L }{\Lambda }
\newcommand{\m }{\mu }
\newcommand{\n }{\nabla }
\newcommand{\s }{\sigma }
\newcommand{\Sig }{\Sigma }
\renewcommand{\t }{\tau }
\newcommand{\var }{\varphi }
\renewcommand{\o }{\omega }
\renewcommand{\O }{\Omega }
\newcommand{\R}{{\mathbb{R}}}
\newcommand{\bC}{{\bf C}}
\newcommand{\bZ}{{\bf Z}}
\newcommand{\bN}{{\bf N}}
\newcommand{\bQ}{{\bf Q}}
\newcommand{\bK}{{\bf K}}
\newcommand{\bI}{{\bf I}}
\newcommand{\bv}{{\bf v}}
\newcommand{\bV}{{\bf V}}
\DeclareMathOperator{\suppo}{supp} \DeclareMathOperator{\di}{div}

\newenvironment{Proof}{\Rmovelastskip\vskip12pt
plus 1pt \noindent\em\rm}{\hfill {\qed \hskip .2cm}}

\begin{document}

\title[]{Uniqueness in weighted $l^p$ spaces \\ for the Schr\"odinger equation \\ on infinite graphs}

\author{Giulia Meglioli}

\address{\hbox{\parbox{5.7in}{\medskip \noindent{Giulia Meglioli, \\Fakult\"at f\"ur Mathematik, \\Universit\"at Bielefeld, \\33501, Bielefeld, Germany \\ [3pt] \emph{E-mail address: }{\tt gmeglioli@math.uni-bielefeld.de}}}}}

\author{Fabio Punzo}

\address{\hbox{\parbox{5.7in}{\medskip \noindent{Fabio Punzo, \\Dipartimento di Matematica, \\Politecnico di Milano, \\Piazza Leonardo da Vinci 32, 20133, Milano, Italy \\ [3pt] \emph{E-mail address: }{\tt fabio.punzo@polimi.it}}}}}


\subjclass[2010]{}

\maketitle

\begin{abstract} We investigate uniqueness of solutions to {\it Schr\"odinger-type} elliptic equations posed on infinite graphs. Solutions are assumed to belong to suitable weighted $\ell^p$ spaces where $p\geq 1$ and the weight is related to both the potential and $p$.

\end{abstract}

\section{Introduction}\label{sec0}

We are concerned with uniqueness of solutions of the following {\it Schr\"odinger-type}  equation
\begin{equation}\label{problema}
\Delta u -V u=0 \quad \text{ in }\; G,
\end{equation}
where $(G, \omega, \mu)$ is an {\it infinite weighted} graph with {\it edge-weight} $\omega$ and {\it node (or vertex) measure} $\mu$, the potential $V$ is a positive function defined in $G$, $\Delta$ denotes the Laplace operator on $G$.

Graphs arises to model many phenomena in applied sciences (see, e.g., \cite{Deo, Les, LHN})\,. Recently, elliptic and parabolic equations posed on graphs have been the object of a detailed investigation. For instance, elliptic equations have been studied in \cite{GLY1, GLY2, GHY, HW, KS, LY, LY2, PinaS}, whereas in \cite{BCG, CGZ, EM, HL, HMu, LW2, Mu, Wu} parabolic equations have been addressed. Moreover, we mention the monographs \cite{Grig2, Mu2, KLW}.

Let us now briefly recall some results in the literature somehow related to our problem.
In \cite{Huang} the parabolic problem
\[\begin{cases}
\partial_t u = \Delta u  & \text{ in } G\times (0, \infty),\\
u = 0 & \text{ in } G\times \{0\}\,
\end{cases} \]
has been investigated. In particular, under suitable assumptions on the graph, it is shown that it admits at most one solution $u$ fulfilling a suitable growth condition.
To be specific, it is assumed that, for some $C>0$, $0<c<\frac 12$ and a diverging sequence $\{R_n\}_{n\in\mathbb{N}}\subset(0,+\infty)$,
\[\int_0^T\sum_{x\in B_{R_n}(x_0)}u^2(x,t) \mu(x)\,dt \leq C e^{c R_n \log R_n} \quad \text{ for any $n$ large enough.}\]
Such result has been generalized in \cite{HKS}, under more general growth conditions on the solution. A similar result has also been stated in \cite{CGZ} for the \textit{discrete} heat equation.
Analogous uniqueness results for parabolic equations are also known on manifolds (see, e.g., \cite{Grig3, Grig, GrigHK, Pu1}) and in bounded domains of $\mathbb R^n$ (see \cite{NP}, \cite{Pu2}), and for nonlocal operators (see \cite{MP, PV}).

Uniqueness of solutions to Schr\"odinger-type equations like \eqref{problema} has been studied on Riemannian manifolds. More precisely, in \cite{Grig3} bounded solutions are dealt with, instead in \cite{MR} solutions are assumed to belong to a certain weighted $L^p$ space. Analogous results have been obtained in \cite{BP1}, in bounded domains of $\mathbb R^n$.
\medskip

Observe that $\ell^p$-Liouville theorems for equation \eqref{problema} with $V\equiv 0$ have been obtained in \cite{HK}; preliminary results can be found in \cite{Mas, HJ}. These results have been recently generalized also to general regular Dirichlet forms in \cite{HKLS}. Note that Liouville type results for equation \eqref{problema}, with a potential $V$, can also be deduced from Schnol's theorem for $p=2$; a discrete version of such result can be found in \cite{KLW} (see also \cite{LT}). Some comments concerning the connection between such results and our new ones can be found in Remarks \ref{oss5} and \ref{oss6}.

In the present paper we can consider both a potential $V$ and $p\neq 2$. To be specific, we aim at proving that $u\equiv 0$ is the only solution to equation \eqref{problema}, whenever $u$ belongs to a certain $\ell^p_{\varphi}(G, \mu)$ space, where $\varphi$ is a weight which tends to $0$ at {\it infinity}. More precisely, on a suitable class of infinite graphs, we assume that
$u\in \ell^p_{\varphi}(G, \mu)$ for some $p\geq 2$ (see Theorem \ref{teo1}). Furthermore, under stronger hypotheses on the graph, we may also consider $p\in [1, 2)$ (see Theorem \ref{teo2}); indeed, such result is established for any $p\geq 1$, however in the case $p\geq 2$, already covered by the preceding theorem, it requires an extra assumption on the graph. In addition, as a particular case, we also obtain  uniqueness in $\ell^\infty$ (see Remark \ref{oss3}).  Let us mentions that also on Riemannian manifolds for such type of results tipically the case $1<p<2$ is more delicate, and it requires stronger assumptions on the underlying manifold (see \cite{Pu1}).
 Let us stress the fact that, in \cite{CGZ,Huang, HKS}, where parabolic equations are addressed, only $p=2$ is taken into account, instead, we are able to consider also the case of $p\neq2$. Note that, in general, in view of our assumptions,  $\ell^p_{\varphi}(G, \mu)\not\subset \ell^q_{\varphi}(G, \mu)$ when $p<q$.

\noindent Our results are based on suitable a priori estimates (see Propositions \ref{lemma1} and \ref{prop1}). They are obtained by multiplying the equation by appropriate test functions and then integrating by parts. The integration by parts is performed just once when $u\in \ell^p_{\varphi}(G,\mu)$ with $p\ge 2$. On the other hand it is made two times when $u\in \ell^p_{\varphi}(G,\mu)$ with $p\ge 1$. In both cases, the test functions involve functions that must satisfy a suitable inequality. Loosely speaking, they can be thought as supersolution to a certain "adjoint equation", which is different in the two range of the parameter $p$ (see \eqref{eq314} and \eqref{e42f}).
In the construction of such supersolutions, it will be apparent the interplay between  $\inf_G V$ (which is supposed to be positive), the weight $\varphi$ and $p$.

The paper is organized as follows. In Section \ref{mf} we introduce the graph setting, in Section \ref{sec2} we state our main uniqueness results. Section \ref{sec3} is devoted to some auxiliary results, concerning appropriate a priori estimates. The main results are proved in Section \ref{sec4}, for $p\geq 2$, and in Section \ref{sec5}, for $p\geq 1$.

\section{Mathematical framework}\label{mf}
\subsection{The graph setting}
A graph $(G, \omega, \mu)$ is a triplet, where
\begin{itemize}
\item $G$ is a countable set of {\it vertices} or {\it nodes};
\item $\omega:G\times G\to [0, +\infty)$ is a function, called {\it edge weight}:
\item $\mu:G\to (0, +\infty)$ is a function, called {\it node (or vertex) measure};
\end{itemize}
where the function $\omega$ fulfills
\begin{itemize}
\item ({\it no loops}) $\omega(x,x)=0 \quad \text{ for all }\, x\in G$ ;
\item ({\it symmetry}) $\omega(x,y)=\omega(y,x)\quad \text{ for all }\, (x,y)\in G\times G$;
\item ({\it finite sum}) $\displaystyle \sum_{y\in G} \omega(x,y)<\infty \quad \text{ for all }\, x\in G\,.$
\end{itemize}
If $G$ is a finite set, a graph is said to be {\it finite}; if $G$ is infinite and countable, we say that the graph is {\it infinite}. We write $x\sim y$ whenever $\omega(x,y)>0$; this means that the couple $(x,y)$ is an {\it edge} of the graph. In this case, we say that $x$ is {\it connected} to $y$ (or $x$ is {\it joint} to $y$, or $x$ is {\it adjacent} to $y$, or $x$ is a {\it neighbor} of $y$).
The vertices $x,y$ are called the {\it endpoints} of the edge $(x,y)$. We specify that we deal with {\it undirected graphs}, that is the edges do not have an orientation; moreover, it is not possible that an edge connects a vertex to itself, so the graphs do not possess {\it loops}, nor for two vertices to be connected by more than one edge. A collection of vertices $\gamma\equiv \{x_k\}_{k=0}^n\subset G$ is called a {\it path} if $x_k\sim x_{k+1}$ for all $k=0, \ldots, n-1.$ A graph is {\it connected} if, for any two vertices $x,y\in V$, there exists a path joining $x$ to $y$. We say that a graph  is {\it locally finite}, if, for each vertex $x\in G$, $\operatorname{card}\{y\in G\,:\, x\sim y\}<\infty\,.$
We define the {\it degree} $\operatorname{deg}$ and the {\it weighted degree} $\operatorname{Deg}$ of a vertex $x\in G$ by
\[\operatorname{deg}(x):=\sum_{y\in G}\omega(x,y), \quad \operatorname{Deg}(x):=\frac{\operatorname{deg}(x)}{\mu(x)}\,.\]

\noindent A {\it pseudo metric} on $G$ is a mapping $d:G\times G\to [0, +\infty)$ such that
\begin{itemize}
  \item $d(x,x)=0$ \quad \text{ for all } $x\in G$;
  \item $d(x,y)=d(y,x)$ \text{ for all } $x,y\in G$;
  \item $d(x,y)\leq d(x,z)+d(z,y)$ \text{ for all } $x,y,z\in G$.
\end{itemize}
In general, $d$ is not a metric, since we can find points $x, y\in G, x\neq y$ such that $d(x,y)=0\,.$
\begin{definition}
We say that
\begin{equation}\label{e14f}
s:=\sup\{d(x,y) \,:\, x,y\in G, \omega(x,y)>0\}
\end{equation}
is the {\em jump size} of $d$.
\end{definition}

\smallskip

\begin{definition}\label{intrmet}
Let $q\geq 1, C_0>0.$ We say that a pseudo metric $d$ on $(G, \omega, \mu)$ is $q-${\em intrinsic} with bound $C_0$ if
\begin{equation}\label{e1f}
\frac 1{\mu(x)}\sum_{y\in G} \omega(x,y)d^q(x,y) \leq C_0\quad \text{ for all } x\in G\,.
\end{equation}
If $q=2, C_0=1$, then $d$ is usually called {\em intrinsic} (see e.g. \cite{HKS}).
\end{definition}


In Remark \ref{oss2} we comment about condition \eqref{e1f}. 

\subsection{Difference and Laplace operators} Let $\mathfrak F$ denote the set of all functions $f: G\to \mathbb R$\,. For any $f\in \mathfrak F$ and for all $x,y\in G$, let us introduce the {\it difference operator}
\begin{equation}\label{e2f}
\nabla_{xy} f:= f(y)-f(x)\,.
\end{equation}
It is direct to see that for any $f,g\in \mathfrak F$, the next {\it product rule} takes place
\begin{equation}\label{e3f}
\nabla_{xy}(fg)=f(x) (\nabla_{xy} g) + (\nabla_{xy} f)g(y) \quad \text{ for all } x,y\in G\,.
\end{equation}
The gradient squared of $f\in \mathfrak F$ is defined as (see \cite{CGZ})
\[|\nabla f(x)|^2=\frac 1{\mu(x)}\sum_{y\in G}\omega(x,y)(\dif f)^2, \quad x\in G\,.\]

\begin{definition}\label{def1}
Let $(G, \omega,\mu)$ be a weighted graph. For any $f\in \mathfrak F$ such that 
\begin{equation}\label{sumfinite}
\sum_{y\in G}\omega(x, y) |f(y)|<\infty \quad \text{ for all } \; x\in G,
\end{equation}
 the {\em (weighted) Laplace operator} of $(G, \omega, \mu)$ is
\begin{equation}\label{eq11}
\Delta f(x):=\frac{1}{\mu(x)}\sum_{y\in G}[f(y)-f(x)]\omega(x,y)\quad \text{ for all }\, x\in G\,.
\end{equation}
\end{definition}
Note that \eqref{sumfinite} holds whenever the graph is locally finite. Clearly,
\[\Delta f(x)=\frac 1{\mu(x)}\sum_{y\in G}(\dif f)\omega(x,y)\quad \text{ for all } x\in G\,.\]

It is easily seen that if at least one of the functions $f, g\in \mathfrak F$ has {\it finite} support, then the following {\it integration by parts formula} holds:
\begin{equation}\label{e4f}
\sum_{x\in G}[\Delta f(x)] g(x) \mu(x)=-\frac 1 2\sum_{x,y\in G}(\dif f)(\dif g) \omega(x,y)\,.
\end{equation}

\begin{remark}\label{oss2}
Now we compare hypothesis \eqref{e1f} with analogous condition on Riemannian manifolds. Let $x_0\in G$ be any fixed reference point, consider the function $x\mapsto d(x, x_0)$ for $x\in G.$  Then
\[|\nabla d(x, x_0)|^2=\frac 1{\mu(x)}\sum_{y\in G}\omega(x,y)[d(y, x_0)-d(x, x_0)]^2\leq \frac 1{\mu(x)} \sum_{y\in G}\omega(x,y) d^2(x,y)\;  \text{ for any } x\in G \,.\]
Therefore, hypothesis \eqref{e1f} (with $q=2, C_0=1)$ ensures $|\nabla d(x, x_0)|^2\leq 1$ for any $x\in G\,.$
Such property is clearly fulfilled on Riemannian manifolds. Moreover,
\[\big|\Delta d(x, x_0) \big|\leq \frac 1{\mu(x)} \sum_{y\in G}\omega(x,y)\big|d(y, x_0)-d(x, x_0)\big|\leq \frac 1{\mu(x)} \sum_{y\in G}\omega(x,y) d(x,y)  \quad \text{ for any }\; x\in G\,. \]
Hence, \eqref{e1f} (with $q=1)$ yields
\begin{equation}\label{e65f}
\big|\Delta d(x, x_0) \big|\leq C_0 \quad \text{ for any }\; x\in G\,.
\end{equation}
On Riemannian manifolds, from Laplacian and Hessian comparison results it follows that \eqref{e65f} is fulfilled, whenever the Ricci curvature is bounded from below and the sectional curvature is bounded from above.
\end{remark}

\subsection{Weighted $\ell^p$ spaces}\label{ellep} Let $\varphi:G\to (0, +\infty)$ be a function. For any $p\in [1, +\infty]$ we define the $\ell^p_{\varphi}(G, \mu)$ space as
\[\ell^p_{\varphi}(G, \mu):=\{u\in \mathfrak F\,:\, \sum_{x\in G}|u(x)|^p \varphi(x) \mu(x) <+\infty\} \quad \text{ for } p\in [1, +\infty),\]
\[\ell^\infty_{\varphi}(G, \mu)\equiv \ell^\infty(G):=\{u\in \mathfrak F\,:\, \sup_{x\in G}|u(x)| <+\infty\}.\]

When $\varphi\equiv 1$, we get the usual $\ell^p(G, \mu)$ spaces (see \cite{BSW}).

\section{Main results} \label{sec2}\setcounter{equation}{0}

\noindent For any $x_0\in G$ and $r>0$ we set
\[B_r(x_0):=\{x\in G\,:\, d(x,x_0)<r\}\,.\]

\noindent We always make the following assumption:
\begin{equation}\label{e7f}
\begin{cases}
(i) & (G, \mu, \omega) \text{ is an infinite, connected, weighted graph};\\
(ii) & \text{there exists a pseudo distance } d \text{ such that the jump size $s$ is finite,  } \\
& \text{and } B_r(x) \text{ is a finite set,  for any } x\in G, r>0.\\
\end{cases}
\end{equation}
Note that by the previous hypothesis, $G$ is automatically locally finite. 

\noindent Concerning the {\it potential} $V$, we always suppose that
\begin{equation}\label{e12f}
V\in \mathfrak F,\;\,\, c_0:= \inf_G V>0.
\end{equation}

Let $x_0\in G$. For each parameter $\gamma>0$, we define the function
\begin{equation}\label{e60f}
\varphi_\gamma(x):=e^{-\gamma d(x, x_0)}\quad \text{ for any }\; x\in G\,.
\end{equation}

\begin{theorem}\label{teo1}
Assume that $d$ is an intrinsic pseudo metric with finite balls and finite jump size $s$. Furthermore, assume that $c_0:= \inf_G V>0$. Let $u$ be a solution to problem \eqref{problema}. Let $p\geq 2$
and $\beta>0$ be such that
\begin{equation}\label{e13f}
\beta^2 e^{2 s \beta}<2 c_0 p\,.
\end{equation}
Suppose that $u\in \ell^p_{\varphi_\beta}(G, \mu).$ Then $ u(x) = 0 \; \text{ for any } \; x\in G.$
\end{theorem}
 \normalcolor

\begin{remark}\label{oss1} (i) From the proof of Theorem \ref{teo1} it is easily seen that instead of $u\in \ell^p_{\varphi_\beta}(G, \mu)$ it could be assumed that
there exist $C>0$ and a sequence $\{R_n\}_{n\in\mathbb{N}}\subset (1, +\infty)$ with
$\displaystyle R_n\mathop{\longrightarrow}_{n\to +\infty}+\infty$, such that
\begin{equation*}
\sum_{x\in B_{R_n}(x_0)}|u(x)|^p\mu(x)\,\le\, C e^{\beta R_n}\quad \text{ for all }  n\in \mathbb N\,.
\end{equation*}

\noindent (ii) In general, observe that if $u\in \ell^p_{\varphi_\beta}(G, \mu)$, then, for some $C>0$,
\begin{equation}\label{h2}
 \sum_{x\in B_{R}(x_0)}|u(x)|^p \mu(x) \leq C e^{\beta R}\quad \text{ for any } R>1\,.
 \end{equation}
In fact, $u\in \ell^p_{\varphi_\beta}(G, \mu)$ if and only if
$\displaystyle \sum_{x\in G}|u(x)|^p \mu(x) e^{-\beta d(x_, x_0)}\leq C,$
for some $C>0$. Hence
\[C\geq \sum_{x\in B_{R}(x_0)}|u(x)|^p \mu(x) e^{-\beta d(x_, x_0)}\geq e^{-\beta R}\sum_{x\in B_{R}(x_0)}|u(x)|^p \mu(x) \,. \]
Therefore, \eqref{h2} follows.
\end{remark}

\smallskip

In order to cover also the case $p\in [1, 2)$, we need an hypothesis on $d$; indeed, we also suppose that $d$ is $1-${\em intrinsic} with bound $C_0$. 
\begin{theorem}\label{teo2}
Assume that $d$ is an intrinsic pseudo metric with finite balls and finite jump size $s$, which is also $1-$intrinsic with bound $C_0$. Let $u$ be a solution to problem \eqref{problema}. Furthermore, assume that $c_0:= \inf_G V>0$. Let $p\geq 1$.
Let $\alpha>0$ be such that
\begin{equation}\label{e13fb}
C_0 \alpha e^{ s \alpha}< c_0 p\,. \quad
\end{equation}
Suppose that $u\in \ell^p_{\varphi_\alpha}(G, \mu)$.Then $u(x) = 0 \, \text{ for any } \; x\in G.$
\end{theorem}

Observe that, due to \eqref{e12f}, condition \eqref{e13f} hilights the relation between the potential $V$ and the uniqueness class $\ell^p_{\varphi_\beta}(G, \mu)$; an analogous comment can be done for \eqref{e13fb}.

\begin{remark}
Note that in Theorem \ref{teo2} we cover not only the case $p\in [1, 2)$, but also the case $p\geq 2$ that is already addressed in Theorem \ref{teo1}. However, in Theorem \ref{teo2} we need the extra hypothesis that $d$ is $1-intrinsic$, not required in Theorem \ref{teo1}.
\end{remark}

\begin{remark}\label{oss3}[Uniqueness in $\ell_\varphi^\infty(G, \mu)$]
Note that, in general, if we assume that
\begin{equation}\label{e66f}
\sum_{x\in G} \varphi_{\gamma} (x)\mu(x)<+\infty,
\end{equation}
then $\ell^\infty_\varphi(G, \mu)\subset \ell^p_\varphi(G, \mu),$ for each $p\geq 1$.  Hence, in Theorem \ref{teo1}, if we assume that \eqref{e66f} holds with $\gamma=\beta$, then it follows uniqueness of solutions to \eqref{problema} in $\ell^\infty(G, \mu)$. The same holds for Theorem \ref{teo2}, if we assume \eqref{e66f} with $\gamma=\alpha$.
\end{remark}

\begin{remark}\label{oss4}
As it can be easily seen from the proofs, in Theorems \ref{teo1} and \ref{teo2}, instead of requiring that $u$ is a solution of equation \eqref{problema}, we can require that $u\geq 0$ and $\Delta u \geq c_0 u$. 
\end{remark}

\begin{remark}\label{oss5}
Note that in \cite{HK} a discrete version of Karp's theorem is established. It is shown that if $u\geq 0$ fulfills $\inf_{r_0>0}\int_{r_0}^{+\infty}\frac r{\|u {\bf 1}_{B_r}\|_p^p}dr=+\infty$ and $\Delta u \geq 0$, then $u\equiv 0$.
Due to Theorems \ref{teo1}-\ref{teo2} and Remark \ref{oss4}, we can say that by replacing $\Delta u\geq 0$ by $\Delta u\geq c_0 u$, we can allow some additional exponential growth on $u$ However, by our methods we are not able to consider an integral condition like that in Karp's theorem.
\end{remark}

\begin{remark}\label{oss6}
From Schnol's theorem \cite[Theorem 12.25]{KLW} it follows that if $u$ solves equation \eqref{problema} and $u\in \ell^2_{\varphi_\beta}(G, \mu)$ for all $\beta>0$, then $u\equiv 0$. Observe that in Theorem \ref{teo1} the request on $u$ is made only for some small enough $\beta>0$. The methods of proof of the two results are totally different.
\end{remark}

\section{Auxiliary Results}\label{sec3}\setcounter{equation}{0}
In all this Section $(G, \mu, \omega)$ is an infinite, connected, locally finite, weighted graph.

The following two lemmas will be expedient in the sequel.

\begin{lemma}\label{lemma4} 
Let $\psi\in C^1(\mathbb R; \mathbb R)$ be a convex function, let $u\in \mathfrak F$. Then
\begin{equation}\label{e5f}
\Delta \psi(u(x)) \geq \psi'(u(x)) \Delta u(x)\quad \text{ for all }\; x\in G\,.
\end{equation}
\end{lemma}
\begin{proof}
Since $\psi\in C^1(\mathbb R;\mathbb R)$ is convex, we have that
\begin{equation}\label{e6f}
\psi(s)- \psi(t) \geq \psi'(t)(s-t)\quad \text{ for any }\; s,t\in\R\,.
\end{equation}
Let $x\in G$. Due to \eqref{eq11} and \eqref{e6f} with $s=u(y)$ and $t=u(x)$, we get
\[
\begin{aligned}
\Delta \psi(u(x))&= \frac1{\mu(x)}\sum_{y\in G}[\psi(u(y)- \psi(u(x))]\omega(x,y)\\&\ge\frac{1}{\mu(x)}\sum_{y\in G} \psi'(u(x))\left[u(y)-u(x)\right]\omega(x,y)\\&
= \frac{\psi'(u(x))}{\mu(x)}\sum_{y\in G} \left[u(y)-u(x)\right]\omega(x,y) =\psi'(u(x))\Delta u(x)\,.
\end{aligned}
\]
\end{proof}

\begin{lemma}[\textit{Laplacian of the product}]\label{lemma5} 
Let $f,g\in \mathfrak F$. Then, for every $x\in G$,
\begin{equation}\label{e25f}
\Delta [f(x)g(x)]=f(x) \Delta g(x) + g(x) \Delta f(x) + \frac 1{\mu(x)}\sum_{y\in G}(\dif f)(\dif g )\omega(x,y)\,.
\end{equation}
\end{lemma}
\begin{proof}
Let $f,g\in \mathfrak F.$ By \eqref{e3f} and \eqref{eq11},
\begin{equation}\label{e26f}
\Delta[f(x) g(x)]=\frac 1{\mu(x)}\sum_{y\in G}(\dif fg) \omega(x,y)=\frac1{\mu(x)}\sum_{y\in G}\left[f(x) (\dif g) + g(y) (\dif f) \right]\omega(x,y)
\end{equation}
for all $x\in G$. Observe that, for every $x\in V,$
\begin{equation}\label{e27f}
\begin{aligned}
\sum_{y\in G}g(y)(\dif f)\omega(x,y)&=\sum_{y\in G}[g(y)-g(x)+g(x)](\dif f)\omega(x,y)\\
&=g(x)\sum_{y\in G}(\dif f)\omega(x,y)+\sum_{y\in G}(\dif g)(\dif f)\omega(x,y)\,.
\end{aligned}
\end{equation}
From \eqref{e26f} and \eqref{e27f} we easily obtain \eqref{e25f}.
\end{proof}

The following proposition is expedient for the proof of Theorem \ref{teo1}.

\begin{proposition}\label{lemma1}
Let $u$ be a solution to equation \eqref{problema}. Let $\eta, \xi\in \mathfrak F$; suppose that
\begin{align}
&\text{$\bullet$}\quad \eta\ge0, \; \operatorname{supp} \eta \text{ is finite };\nonumber
\\
&\text{$\bullet$}\quad [\eta^2(y)-\eta^2(x)][e^{\xi(y)}-e^{\xi(x)}]\ge0 \quad \text{ for all } x,y\in G, x\sim y.\label{e10f}
\end{align}
Then, for any $p\ge 2$,
\begin{equation*}\label{eq31}
\begin{aligned}
\frac 12 \sum_{x\in G}|u(x)|^p\eta^2(x)&e^{\xi(x)}\left\{V(x) \,p \mu(x)-\frac 12\sum_{y\in G}\omega(x,y)\left[1-e^{\xi(y)-\xi(x)}\right]^2 \right\}\\
&\le\sum_{x,y\in G}|u(x)|^pe^{\xi(y)}\left[\eta(y)-\eta(x)\right]^2\omega(x,y)\,.
\end{aligned}
\end{equation*}
\end{proposition}

\begin{proof}
For any $\alpha>0$, we define
\begin{equation}\label{G}
\pi_{\alpha}(u):=(u^2+\alpha)^{\frac p4}\,.
\end{equation}
Now, for any $x\in G$, we consider the expression
\begin{equation}\label{eq32}
\Delta \pi_{\alpha}(u(x))-V(x)\,\pi_{\alpha}(u(x)).
\end{equation}
We multiply \eqref{eq32} by $\pi_{\alpha}(u(x))\,\eta^2(x)e^{\xi(x)}\mu(x)$ and then we sum over $x\in G$. Thus
\begin{equation}\label{eq33}
\sum_{x\in G}\Delta \pi_{\alpha}(u(x))\pi_{\alpha}(u(x))\,\eta^2(x)e^{\xi(x)}\mu(x)\,-\sum_{x\in G}\,V(x)\,\pi^2_{\alpha}(u(x))\,\eta^2(x)e^{\xi(x)}\mu(x)\,.
\end{equation}
Note that since $\eta(x)$ is finitely supported, the series in \eqref{eq33} are finite sums. Set
\begin{equation}\label{eq34a}
I:=\sum_{x\in V}\Delta \pi_{\alpha}(u(x))\pi_{\alpha}(u(x))\,\eta^2(x)e^{\xi(x)}\mu(x)\,.
\end{equation}
Then, by \eqref{e4f} and \eqref{e2f},
\begin{equation}\label{eq34}
\begin{aligned}
I&=-\frac 12 \sum_{x,y\in G}\left(\nabla_{xy} \pi_{\alpha}(u)\right) \nabla_{xy}\left[\pi_{\alpha}(u)\eta^2e^{\xi}\right]\omega(x,y)\\
&=-\frac 12 \sum_{x,y\in G}\left[\nabla_{xy}\pi_{\alpha}(u)\right]^2\eta^2(y)e^{\xi(y)}\,\omega(x,y)\\
&\,\,\,\,\,-\frac 12 \sum_{x,y\in G}\pi_{\alpha}(u(x))\eta^2(x)\left(\nabla_{xy}\pi_{\alpha}(u)\right)\left(\nabla_{xy}e^{\xi}\right)\,\omega(x,y)\\
&\,\,\,\,\,-\frac 12 \sum_{x,y\in G}\pi_{\alpha}(u(x))e^{\xi(y)}\left(\nabla_{xy}\pi_{\alpha}(u)\right)\left(\nabla_{xy}\eta^2\right)\,\omega(x,y)\\
&=:J_1+J_2+J_3\,.
\end{aligned}
\end{equation}
In view of \eqref{e2f}, we obviously have
\begin{equation}\label{e11f}
\dif \eta^2=[\eta(y)+\eta(x)][\eta(y)-\eta(x)] \quad \text{ for all }\, x,y\in G\,.
\end{equation}
Due to \eqref{e11f}, by Young's inequality with exponent $2$, we have, for any $\delta_1>0$,
\begin{equation}\label{eq35}
\begin{aligned}
J_3&= -\frac 12 \sum_{x,y\in G}\pi_{\alpha}(u(x))e^{\xi(y)}(\nabla_{xy}\pi_{\alpha}(u))\left[\eta(y)+\eta(x)\right]\left[\eta(y)-\eta(x)\right]\,\omega(x,y)\\
&\le \frac{\delta_1}{4} \sum_{x,y\in G}e^{\xi(y)}\left[\nabla_{xy}\pi_{\alpha}(u)\right]^2\left[\eta(y)+\eta(x)\right]^2\,\omega(x,y)\\
&\,\,\,\,\,+\frac 1{4\delta_1} \sum_{x,y\in G}e^{\xi(y)} \pi_{\alpha}^2(u(x))\left[\eta(y)-\eta(x)\right]^2\,\omega(x,y)\\
&\le \frac{\delta_1}{2} \sum_{x,y\in G}e^{\xi(y)}\left[\nabla_{xy}\pi_{\alpha}(u)\right]^2 \left[\eta^2(y)+\eta^2(x)\right]\,\omega(x,y)\\
&\,\,\,\,\,+\frac 1{4\delta_1} \sum_{x,y\in G}e^{\xi(y)} \pi_{\alpha}^2(u(x))\left[\eta(y)-\eta(x)\right]^2\,\omega(x,y)\\
&= \frac{\delta_1}{4} \sum_{x,y\in G}\left[e^{\xi(y)}+e^{\xi(x)}\right]\left[\nabla_{xy}\pi_{\alpha}(u)\right]^2\left[\eta^2(y)+\eta^2(x)\right]\,\omega(x,y)\\
&\,\,\,\,\,+\frac 1{4\delta_1} \sum_{x,y\in G}e^{\xi(y)} \pi_{\alpha}^2(u(x))\left[\eta(y)-\eta(x)\right]^2\,\omega(x,y)\,.\\
\end{aligned}
\end{equation}
From \eqref{e10f} we can easily infer that
\begin{equation}\label{eq36}
\left[\eta^2(y)+\eta^2(x)\right]\left[e^{\xi(y)}+e^{\xi(x)}\right]\le 2\left\{\eta^2(y)e^{\xi(y)}+\eta^2(x)e^{\xi(x)}\right\}\,.
\end{equation}
By using \eqref{eq35} and \eqref{eq36},
\begin{equation}\label{eq37}
\begin{aligned}
J_3&\le \frac{\delta_1}{2} \sum_{x,y\in G}\left[\nabla_{xy}\pi_{\alpha}(u)\right]^2\left\{\eta^2(y)e^{\xi(y)}+\eta^2(x)e^{\xi(x)}\right\}\,\omega(x,y)\\
&\,\,\,\,\,+\frac 1{4\delta_1} \sum_{x,y\in G}e^{\xi(y)} \pi_{\alpha}^2(u(x))\left[\eta(y)-\eta(x)\right]^2\,\omega(x,y)\\
&=\delta_1\sum_{x,y\in G}\left[\nabla_{xy}\pi_{\alpha}(u)\right]^2\eta^2(y)e^{\xi(y)}\,\omega(x,y)\\
&\,\,\,\,\,+\frac 1{4\delta_1} \sum_{x,y\in G}e^{\xi(y)} \pi_{\alpha}^2(u(x))\left[\eta(y)-\eta(x)\right]^2\,\omega(x,y)\,.\\
\end{aligned}
\end{equation}
Similarly, by Young's inequality with exponent $2$, for any $\delta_2>0$, we have
\begin{equation}\label{eq38}
\begin{aligned}
J_2&=-\frac 12 \sum_{x,y\in G}\pi_{\alpha}(u(x))\eta^2(x)(\nabla_{xy}\pi_{\alpha}(u))(\nabla_{xy}e^{\xi})\,\omega(x,y)\\
&=-\frac 12 \sum_{x,y\in G}\pi_{\alpha}(u(y))\eta^2(y)(\nabla_{xy}\pi_{\alpha}(u))\left[e^{\xi(y)}-e^{\xi(x)}\right]\,\omega(x,y)\\
&=-\frac 12 \sum_{x,y\in G}\pi_{\alpha}(u(y))\eta^2(y)\nabla_{xy}\pi_{\alpha}(u)e^{\xi(y)}\left[1-e^{\xi(x)-\xi(y)}\right]\,\omega(x,y)\\
&\le \frac{\delta_2}{4}\sum_{x,y\in G}\left[(\nabla_{xy}\pi_{\alpha}(u))\right]^2\eta^2(y)e^{\xi(y)}\,\omega(x,y)\\
&\quad + \frac1{4\delta_2}\sum_{x,y\in G}\pi^2_{\alpha}(u(y))\left[1-e^{\xi(x)-\xi(y)}\right]^2\eta^2(y)e^{\xi(y)}\,\omega(x,y)\,.
\end{aligned}
\end{equation}
By combining together \eqref{eq34}, \eqref{eq37} and \eqref{eq38} we get
\begin{equation}\label{eq39}
\begin{aligned}
\sum_{x\in G}&\Delta \pi_{\alpha}(u(x))\pi_{\alpha}(u(x))\,\eta^2(x)e^{\xi(x)}\mu(x)\\
&\le -\frac 12 \sum_{x,y\in G}[\nabla_{xy}\pi_{\alpha}(u)]^2\eta^2(y)e^{\xi(y)}\omega(x,y) + \frac{\delta_2}{4}\sum_{x,y\in G}[\nabla_{xy}\pi_{\alpha}(u)]^2\eta^2(y)e^{\xi(y)}\omega(x,y)\\
&\quad + \frac1{4\delta_2}\sum_{x,y\in G}\pi^2_{\alpha}(u(y))\left[1-e^{\xi(x)-\xi(y)}\right]^2\eta^2(y)e^{\xi(y)}\,\omega(x,y) \\
&\quad + \delta_1\sum_{x,y\in G}\left[\nabla_{xy}\pi_{\alpha}(u)\right]^2\eta^2(y)e^{\xi(y)}\,\omega(x,y)\\
&\quad +\frac 1{4\delta_1} \sum_{x,y\in G} \pi_{\alpha}^2(u(x))\left[\eta(y)-\eta(x)\right]^2e^{\xi(y)}\,\omega(x,y)\,.
\end{aligned}
\end{equation}
We now choose $\delta_1=\frac 14$ and $\delta_2=1$ in such a way that $-\frac 12+ \frac{\delta_2}{4}+ \delta_1=0$. Consequently, \eqref{eq39} yields
\begin{equation}\label{eq310}
\begin{aligned}
\sum_{x\in G}\Delta \pi_{\alpha}&(u(x))\pi_{\alpha}(u(x))\,\eta^2(x)e^{\xi(x)}\mu(x)\\
&\le \frac1{4}\sum_{x,y\in G}\pi^2_{\alpha}(u(x))\left[1-e^{\xi(y)-\xi(x)}\right]^2 \eta^2(x)e^{\xi(x)}\,\omega(x,y) \\
&\quad + \sum_{x,y\in G}\pi_{\alpha}^2(u(x))\left[\eta(y)-\eta(x)\right]^2e^{\xi(y)}\,\omega(x,y)\,.
\end{aligned}
\end{equation}
By substituting \eqref{eq310} into \eqref{eq33}, we get
\begin{equation}\label{eq311}
\begin{aligned}
\sum_{x\in G}&\Delta \pi_{\alpha}(u(x))\pi_{\alpha}(u(x))\,\eta^2(x)e^{\xi(x)}\mu(x)-\sum_{x\in G}\,V(x)\,\pi^2_{\alpha}(u(x))\,\eta^2(x)e^{\xi(x)}\mu(x)\\
&\le \frac1{4}\sum_{x,y\in G}\pi^2_{\alpha}(u(x))\left[1-e^{\xi(y)-\xi(x)}\right]^2 \eta^2(x)e^{\xi(x)}\,\omega(x,y) \\
&\quad + \sum_{x,y\in G}\pi_{\alpha}^2(u(x))[\eta(y)-\eta(x)]^2e^{\xi(y)}\omega(x,y) -\sum_{x\in G}V(x)\pi^2_{\alpha}(u(x))\,\eta^2(x)e^{\xi(x)}\mu(x)\,.
\end{aligned}
\end{equation}
Observe that, due to the assumption $p\ge2$, $s\mapsto \pi_{\alpha}(s)$ defined in \eqref{G} is convex in $\R$ and of class $C^2$. Therefore, from \eqref{e5f} we get
\begin{equation}\label{eq314b}
\Delta \pi_{\alpha}(u(x))\ge \pi_{\alpha}'(u(x))\Delta u(x)\quad \text{ for all }\; x\in G\,.
\end{equation}
Thus, due to \eqref{eq34a}, \eqref{eq314b} and \eqref{problema}, we can infer that
\begin{equation}\label{eq316b}
\begin{aligned}
I&\ge \sum_{x\in G}\pi_{\alpha}(u(x))\pi'_\alpha(u(x))\Delta u(x)\eta^2(x)e^{\xi(x)}\mu(x)\\
&= \sum_{x\in G}\pi_{\alpha}(u(x))\pi'_\alpha(u(x))V(x) u(x)\eta^2(x)e^{\xi(x)}\mu(x)\,.
\end{aligned}
\end{equation}
By combining together \eqref{eq311} and \eqref{eq316b} we get
\begin{equation*}
\begin{aligned}
\sum_{x\in G}\pi_{\alpha}(u(x))&V(x)\eta^2(x)e^{\xi(x)}\mu(x)\left[\pi'_\alpha(u(x))u(x)-\pi_{\alpha}(u(x))\right]\\
&\le \frac1{4}\sum_{x,y\in G}\pi^2_{\alpha}(u(x))\left[1-e^{\xi(y)-\xi(x)}\right]^2 \eta^2(x)e^{\xi(x)}\,\omega(x,y) \\
&\quad + \sum_{x,y\in G}\pi_{\alpha}^2(u(x))\left[\eta(y)-\eta(x)\right]^2e^{\xi(y)}\,\omega(x,y) -\sum_{x\in G}\,V(x)\,\pi^2_{\alpha}(u(x))\,\eta^2(x)e^{\xi(x)}\mu(x)\,;
\end{aligned}
\end{equation*}
therefore,
\begin{equation}\label{eq317b}
\begin{aligned}
\sum_{x\in G}\pi_{\alpha}(u(x))&V(x)\eta^2(x)e^{\xi(x)}\mu(x)\left[\pi'_\alpha(u(x))u(x)\right]\\
&\le \frac1{4}\sum_{x,y\in G}\pi^2_{\alpha}(u(x))\left[1-e^{\xi(y)-\xi(x)}\right]^2 \eta^2(x)e^{\xi(x)}\,\omega(x,y) \\
&\quad + \sum_{x,y\in G}\pi_{\alpha}^2(u(x))\left[\eta(y)-\eta(x)\right]^2e^{\xi(y)}\,\omega(x,y)\,.
\end{aligned}
\end{equation}
Now, since
\[\pi_\alpha(u)\to |u|^{\frac p2}, \quad  \pi_{\alpha}'(u)u\to\frac p2 |u|^{\frac p2} \quad \text{ as } \alpha\to 0^+,\]
letting $\alpha\to 0^+$ in \eqref{eq317b}, we arrive to
\begin{equation*}\label{eq312}
\begin{aligned}
\frac p2\sum_{x\in G}|u(x)|^p\eta^2(x)e^{\xi(x)} V(x)\mu(x)&\leq \frac1{4}\sum_{x, y\in G}\left[1-e^{\xi(y)-\xi(x)}\right]^2 \,\omega(x,y)\\
&  + \sum_{x,y\in G}|u(x)|^p\left[\eta(y)-\eta(x)\right]^2e^{\xi(y)}\,\omega(x,y)\,.
\end{aligned}
\end{equation*}
This immediately implies the thesis.
\end{proof}

The following two propositions will be used in the proof of Theorem \ref{teo2}.

\begin{proposition}\label{prop1}
 Let $u$ be a solution of equation \eqref{problema}. Let $p\geq 1$, and $v\in \mathfrak F$ be a nonnegative function with finite support. Then
\begin{equation}\label{e30f}
\sum_{x\in G}|u(x)|^p\left\{-\Delta v(x) + p V(x) v(x) \right\}\mu(x) \leq 0\,.
\end{equation}
\end{proposition}
\begin{proof}
For every $\alpha>0$, define
\begin{equation*}\label{e31f}
\phi_\alpha(u):=(u^2+\alpha)^{\frac p2}\,.
\end{equation*}
By applying twice \eqref{e4f}, we obtain
\begin{equation}\label{e32f}
\sum_{x\in G}\left[\Delta\phi_\alpha(u(x)) - V(x)\phi_\alpha(u(x))\right]v(x) \mu(x)=\sum_{x\in G}\left[\Delta v(x) - V(x)v(x)\right]\phi_\alpha(u(x)) \mu(x)\,.
\end{equation}
Since $p\geq 1$, the function $s\mapsto \phi_\alpha(s)$ is convex in $(0, +\infty)$ and of class $C^2$. Therefore, from \eqref{e5f} we get
\begin{equation}\label{e33f}
\Delta\phi_\alpha(u(x))\geq \phi'_\alpha(u(x))\Delta u(x)\quad \text{ for all }\; x\in G\,.
\end{equation}
Due to \eqref{e32f} and \eqref{e33f}, we obtain
\begin{equation}\label{e34f}
\sum_{x\in G}\{\phi'_\alpha(u(x))\Delta u(x)- V(x) \phi_\alpha(u(x))\}v(x)\mu(x)\leq \sum_{x\in G}[\Delta v(x) - V(x)v(x)]\phi_\alpha(u(x)) \mu(x).
\end{equation}
In view of \eqref{problema}, we have that, for any $x\in V$,
\begin{equation}\label{e35f}
\begin{aligned}
\phi'_\alpha(u(x))\Delta u(x)- V(x) \phi_\alpha(u(x)) &=\phi'_\alpha(u(x)) V(x)\,u(x)- V(x) \phi_\alpha(u(x))\\
&=p u (u^2+\alpha)^{\frac p2 -1} V(x)\, u(x)  - V(x)(u^2+\alpha)^{\frac p2}\\
&=V(x) (u^2+\alpha)^{\frac p2-1}\left[(p-1)u^2-\alpha \right]\,.
\end{aligned}
\end{equation}
From \eqref{e25f}, \eqref{e34f} and \eqref{e35f} it follows that
\begin{equation}\label{e36f}
\sum_{x\in G} v(x) (u^2+\alpha)^{\frac p2-1}V(x)[(p-1)u^2-\alpha]\mu(x)\leq \sum_{x\in G}\phi_\alpha(u(x))[\Delta v(x) - V(x) v(x)]\mu(x).
\end{equation}
Letting $\alpha\to 0^+$ in \eqref{e36f}, we get
\[(p-1)\sum_{x\in G} v(x) |u(x)|^p V(x)\mu(x)\leq \sum_{x\in G}|u(x)|^p\left[\Delta v(x) - V(x) v(x)\right]\mu(x)\,.\]
This easily yields \eqref{e30f}.
\end{proof}

\section{Proof of Theorem \ref{teo1}}\label{sec4}
\setcounter{equation}{0}
In this Section, we always assume the hypotheses of Theorem \ref{teo1}. Let $x_0\in G$, $\alpha>0$, $R>0$, $0<\delta<1$. We define the function
\begin{equation}\label{eq313}
\xi(x)=-\alpha\left[d(x,x_0)-\delta R\right]_+\quad \text{ for all }\; x\in G\,.
\end{equation}

\begin{lemma}\label{lemma2}
Let $p\geq 2$. Let $\xi$ be the function defined in \eqref{eq313}.
Then
\begin{equation}\label{eq314}
\frac 12 \sum_{y\in G} \omega(x,y)\left[1-e^{\xi(y)-\xi(x)}\right]^2 -p V(x)\,\mu(x)\leq \mathcal H \mu(x) \quad \text{ for all }\, x\in G\,,
\end{equation}
where
\begin{equation}\label{e16f}
\mathcal H:=\frac{\alpha^2}2e^{2s \alpha}-c_0 p \,.
\end{equation}
\end{lemma}

\begin{proof}
It is direct to see that
\begin{equation}\label{e58f}
(e^a-1)^2\le e^{2|a|}a^2\quad \text{ for every }\; a\in \mathbb R\,.
\end{equation}
Moreover, by the triangle inequality for the distance $d$,
\[|\xi(x)-\xi(y)|\leq \alpha d(x,y) \quad \text{ for all } x,y\in G\,.\]
Therefore
\begin{equation*}\label{eq315}
\begin{aligned}
\frac 12 \sum_{y\in G} \omega(x,y)\left[1-e^{\xi(y)-\xi(x)}\right]^2  &\le \frac 12 \sum_{y\in G} \omega(x,y)\left[\xi(x)-\xi(y)\right]^2e^{2|\xi(x)-\xi(y)|}\\
&\le \frac 12 \sum_{y\in G}\omega(x,y)\alpha^2d^2(x,y)e^{2\alpha d(x,y)}\\
&\le \frac 12 \alpha^2\mu(x)e^{2s \alpha}\,,
\end{aligned}
\end{equation*}
where we have used \eqref{e14f} and the fact that
$
w(x,y)>0\iff x\sim y\,.
$
Then, thanks to \eqref{e16f}, one has
$$
\frac 12 \sum_{y\in G} \omega(x,y)\left[1-e^{\xi(y)-\xi(x)}\right]^2 -V(x)p\,\mu(x)\le \frac 12 \alpha^2\mu(x)e^{2\alpha}- c_0p\,\mu(x)= \mathcal H \mu(x) \quad \text{ for any }\, x\in G\,.
$$
\end{proof}

Let $x_0\in G$ be fixed. For any $R>0$ and $0<\delta<1$ we define the "cut-off" function
\begin{equation}\label{eq316}
\eta(x):=\min\left\{\frac{\left[R-s-d(x_0,x)\right]_+}{\delta R},\,\,1\right\}\, \quad \text{ for any }\; x\in G\,.
\end{equation}

By minor changes in the proof of \cite[Lemma 2.3]{Huang}, we can show the following
\begin{lemma}\label{lemma3}
The function $\eta$ defined in \eqref{eq316} satisfies
\begin{equation}\label{e55f}
|\dif \eta| \leq \frac 1{\delta R} d(x,y)\chi_{\{(1-\delta)R-2s\le d(x,x_0)\le R\}} \quad \text{ for any }\; x, y\in G, \,\,\omega(x,y)>0\,.
\end{equation}
and
\begin{equation*}\label{eq317}
\sum_{y\in G} \left(\dif \eta\right)^2 \omega(x,y)\le \frac{1}{\delta^2 R^2}\mu(x)\chi_{\{(1-\delta)R-2s\le d(x,x_0)\le R\}} \quad \text{ for any }\; x\in G\,.
\end{equation*}
Furthermore,
\begin{equation}\label{e50f}
|\Delta \eta(x)|\leq \frac {C_0}{\delta R} \chi_{\{(1-\delta)R-2s\le d(x,x_0)\le R\}} \quad \text{ for any }\; x\in G\,.
\end{equation}
provided that $d$ is also $1-$ intrinsic with bound $C_0$.
\end{lemma}

\begin{proof}[Proof of Theorem \ref{teo1}]
In view of assumption \eqref{e13f}, we can find $\alpha$ such that
\begin{equation}\label{h3b}
\alpha^2e^{s2\alpha}<2c_0 p
\end{equation}
and
\begin{equation}\label{h3}
\alpha>\beta\,.
\end{equation}
Then, thanks to \eqref{h3}, we select
 $\delta$ so that
\begin{equation}\label{e22f}
0<\delta< \frac 1 2-\frac{\beta}{2\alpha}\,;
\end{equation}
Let $x_0\in G$ be any fixed point.  Moreover, let $\xi=\xi(x)$ and $\eta=\eta(x)\; (x\in G)$ be defined as in Lemmas \ref{lemma2} and \ref{lemma3}, respectively (see \eqref{eq313} and \eqref{eq316}, respectively), with
\begin{equation}\label{e17f}
R>\max\left\{\frac {2 s}{1-2\delta}, 1 \right\}\,.
\end{equation}
One can easily check that $\xi$ and $\eta$ satisfy the assumptions of Proposition \ref{lemma1}. Therefore, due to Proposition \ref{lemma1}, we get
\begin{equation}\label{eq42}
\begin{aligned}
\sum_{x\in G}|u(x)|^p\eta^2(x)&e^{\xi(x)}\left\{V(x)\,p\mu(x)-\frac 1{2}\sum_{y\in G}\omega(x,y)\left[1-e^{\xi(y)-\xi(x)}\right]^2 \right\}\\
&\le2\sum_{x,y\in G}\omega(x,y)|u(x)|^pe^{\xi(y)}\left[\eta(y)-\eta(x)\right]^2\,.
\end{aligned}
\end{equation}
Due to \eqref{e17f}, we have that, if $d(x,x_0)\le\delta R$, then
\begin{equation}\label{e19f}
\eta(x)= 1\quad\quad\text{and}\quad\quad \xi(x)= 0.
\end{equation}
Moreover, in view of \eqref{e16f}, \eqref{eq314} and \eqref{h3b}, the quantity in the curly brackets of the l.h.s. of \eqref{eq42} fulfills
\begin{equation}\label{e18f}
V(x)\,p-\frac 1{2\mu(x)}\sum_{y\in G}\omega(x,y)\left[1-e^{\xi(y)-\xi(x)}\right]^2 \geq |\mathcal H| >0 \quad \text{for every }\; x\in G\,.
\end{equation}
Concerning the l.h.s. of \eqref{eq42}, from \eqref{e18f} and \eqref{e19f} it follows that
\begin{equation}\label{eq43}
\begin{aligned}
\sum_{x\in G}|u(x)|^p\eta^2(x)e^{\xi(x)}&\left\{V(x)\,p\mu(x)-\frac 1{2}\sum_{y\in G}\omega(x,y)\left[1-e^{\xi(y)-\xi(x)}\right]^2 \right\}\\ \geq & |\mathcal H| \sum_{x\in B_{\delta R}(x_0)}|u(x)|^p\mu(x)\,.
\end{aligned}
\end{equation}

Now, we aim at estimating from above the r.h.s. of \eqref{eq42}\,. To this end, recall that, in view of \eqref{e7f} and \eqref{e14f}, $\omega(x,y)>0$ if and only if $d(x,y)\le s$.
Therefore, we get
\begin{equation}\label{e20f}
\xi(y)\leq \alpha s+\xi(x)\quad \quad \text{ for all } x,y \in G \text{ with } \omega(x,y)>0\,.
\end{equation}
From \eqref{e20f} we obtain
\begin{equation}\label{eq44}
2\sum_{x,y\in G}\omega(x,y)|u(x)|^pe^{\xi(y)}[\eta(y)-\eta(x)]^2\le2e^{s\alpha}\sum_{x,y\in G}|u(x)|^pe^{\xi(x)}\omega(x,y)[\eta(y)-\eta(x)]^2.
\end{equation}
By Lemma \ref{lemma3} and \eqref{eq44},
\begin{equation}\label{eq45}
\begin{aligned}
2\sum_{x,y\in G}&\omega(x,y)|u(x)|^pe^{\xi(y)}\left[\eta(y)-\eta(x)\right]^2\\
&\le\frac{2}{(\delta R)^2}e^{s\alpha}\sum_{x\in G}\mu(x)\chi_{\left\{(1-\delta)R-2s\normalcolor\le d(x,x_0)\le R\right\}}|u(x)|^pe^{-\alpha[d(x,x_0)-\delta R]_+}\\
&\le\frac{2}{(\delta R)^2}e^{3s\alpha-(1-2\delta)\alpha R}\sum_{x\in B_{ R}(x_0)}|u(x)|^p\mu(x)\,. \\ 
\end{aligned}
\end{equation}
From \eqref{eq42},\eqref{eq43} and \eqref{eq45} we get
\begin{equation}\label{eq46}
|\mathcal H| \sum_{x\in B_{\delta R }(x_0)}|u(x)|^p\mu(x)\le\frac{2}{(\delta R)^2}e^{3s\alpha-(1-2\delta)\alpha R}\sum_{x\in B_{R }(x_0)}|u(x)|^p\mu(x)\,.\end{equation}
Recall that since $u\in \ell^p_{\varphi_\beta}(G, \mu)$, \eqref{h2} holds (see Remark \ref{oss1}-(ii)).
By \eqref{h2}, \eqref{eq46} reads
\begin{equation}\label{eq46b}
|\mathcal H| \sum_{x\in B_{R}(x_0)}|u(x)|^p\mu(x)\le\frac{2}{\delta^2 R^{2}}e^{3s \alpha}e^{\left[\beta-(1-2\delta)\alpha\right] R}\,.
\end{equation}
In view of \eqref{e22f},
\begin{equation}\label{e23f}
\beta-(1-2\delta)\alpha<0\,.
\end{equation}
By letting $R\to\infty$ in \eqref{eq46b}, due to \eqref{e23f}, we get
\begin{equation}\label{e21f}
|\mathcal H| \sum_{x\in G}|u(x)|^p\mu(x)\le 0\,.
\end{equation}
Since $\mathcal H>0$ and $\mu(x)>0$ for any $x\in G$, \eqref{e21f} implies that
$$u(x) = 0 \quad \text{ for all } \; x\in G. $$
\end{proof}

\section{Proof of Theorem \ref{teo2}}\label{sec5}\setcounter{equation}{0}
In this Section all the hypotheses of Theorem \ref{teo2} are made. For any fixed $x_0\in G$ and $\alpha>0$, define
\begin{equation}\label{e45f}
\zeta(x):=e^{-\alpha d(x, x_0)}\quad \text{ for any }\; x\in G\,.
\end{equation}

\begin{lemma}\label{lemma6}
Let $p\geq 1$.
Then
\begin{equation}\label{e42f}
 \Delta \zeta(x) - p V(x) \zeta(x) \leq \mathcal K \zeta(x) \quad \text{ for all }\; x\in G\,,
\end{equation}
where $\mathcal K:=C_0 \alpha e^{\alpha s}- p c_0\,,$ where $C_0$ is the bound of the $1$- intrinsic metric $d$. 
\end{lemma}
\begin{proof}
In view of \eqref{eq11}, we have
\begin{equation}\label{e43f}
\begin{aligned}
\Delta \zeta(x)&=\frac1{\mu(x)}\sum_{y\in G}\omega(x,y)\left[e^{-\alpha d(y, x_0)}-e^{-\alpha d(x, x_0)} \right]\\
&=\frac{e^{-\alpha d(x, x_0)}}{\mu(x)} \sum_{y\in G}\omega(x,y)\left[e^{-\alpha d(y, x_0)+\alpha d(x, x_0)} -1\right]\\
&\leq \frac{e^{-\alpha d(x, x_0)}}{\mu(x)} \sum_{y\in G}\omega(x,y)\left[e^{\alpha d(x, y)} -1\right]\quad \text{ for any }\; x\in G\,.
\end{aligned}
\end{equation}
It is direct to see that
\begin{equation}\label{e44f}
e^{t}-1\leq te^t \quad \text{ for all }\; t\in \mathbb R\,.
\end{equation}
From \eqref{e45f}, \eqref{e43f} and \eqref{e44f}, we get
\begin{equation}\label{eq45g}
\begin{aligned}
\Delta \zeta(x) -  p V(x) \zeta(x) &\leq \frac{e^{-\alpha d(x, x_0)}}{\mu(x)} \sum_{y\in G}\omega(x,y)\alpha d(x,y) e^{\alpha d(x,y)} -  p V(x) \zeta(x)  \,.
\end{aligned}
\end{equation}
Finally, by using \eqref{e14f}, \eqref{e12f} and \eqref{e1f}, from \eqref{eq45g} we deduce that
\begin{equation*}
\begin{aligned}
\Delta \zeta(x) -  p V(x) \zeta(x)&\le \zeta(x) \left\{\frac{\alpha}{\mu(x)}e^{\alpha s} \sum_{y\in G}\omega(x,y) d(x,y) - p c_0 \right\}\\
&\leq \zeta(x)\left\{ \frac{\alpha e^{\alpha s}}{\mu(x)} C_0\mu(x) -p c_0\right\} = \mathcal K \zeta(x)\quad \text{ for any }\, x\in G\,.
\end{aligned}
\end{equation*}
\end{proof}

\begin{proof}[Proof of Theorem \ref{teo2}]
Let $\eta\equiv \eta_R$ and $\zeta$ be defined by  \eqref{eq316} and \eqref{e45f}, respectively. From Proposition \ref{prop1} with
$$v(x):=\eta_R(x) \zeta(x) \quad \text{ for all\, } x\in G,$$
and Lemma \ref{lemma5} we can infer that
\begin{equation}\label{e47f}
\sum_{x\in G}|u(x)|^p\Big\{-\Delta \zeta(x) + p V(x) \zeta(x) \Big\}\eta_R(x)\mu(x)\leq I_R^{(1)} + I_R^{(2)}\,,
\end{equation}
where
\[I_R^{(1)}:=\sum_{x\in G}|u(x)|^p \zeta(x) \Delta\eta_R(x)\mu(x), \quad I_R^{(2)}:=\sum_{x, y\in G}(\dif \eta_R)(\dif \zeta)|u(x)|^p \omega(x,y).\]
We claim that
\begin{equation}\label{e48f}
\lim_{R\to +\infty}  I_R^{(1)} =0,
\end{equation}
and
\begin{equation}\label{e49f}
\lim_{R\to +\infty}  I_R^{(2)} =0\,.
\end{equation}
From \eqref{e48f} and \eqref{e49f} we can obtain the conclusion. In fact,
letting $R\to +\infty$ in \eqref{e47f}, we have
\begin{equation}\label{e51f}
\sum_{x\in G}|u(x)|^p\Big\{-\Delta \zeta(x) + p V(x) \zeta(x) \Big\}\mu(x)\leq 0\,.
\end{equation}
 In view of \eqref{e13fb} and \eqref{e42f},
 \begin{equation}\label{e52f}
 -\Delta \zeta(x) + p V(x) \zeta(x)  >0 \quad \text{ for all }\; x\in G	\,.
 \end{equation}
Combining \eqref{e52f} and \eqref{e51f}, we deduce that
$u(x)=0$ for all $x\in G.$ It remains to show \eqref{e48f} and \eqref{e49f}. To this purpose, observe that, by exploiting the definition of $\eta_R$ given in \eqref{eq316}, \eqref{e50f} and assumption \eqref{e1f}, we get
\begin{equation}\label{eq56g}
\begin{aligned}
\big| I_R^{(1)} \big| &=\left|\sum_{x\in G}|u(x)|^p \zeta(x)\Delta \eta_R(x)\mu(x)\right|\\
&\le \sum_{x\in G}|u(x)|^p \zeta(x)\mu(x) |\Delta \eta_R(x)|
\le \frac{C_0}{\delta R}\sum_{x\in G}|u(x)|^p \zeta(x)\mu(x)\,.
\end{aligned}
\end{equation}
Since $u\in \ell^p_{\varphi_\alpha}(G, \mu)$, we deduce that, for some $\overline C>0$,
\begin{equation}\label{h2b}
\sum_{x\in G} |u(x)|^p e^{-\alpha d(x, x_0)} \mu(x) \leq \overline C\,.
\end{equation}
Now, by \eqref{h2b}, \eqref{eq56g} yields
\begin{equation*}\label{e56f}
\begin{aligned}
\big| I_R^{(1)} \big| \leq \frac{C_0\overline C}{\delta R} \mathop{\longrightarrow}_{R\to +\infty} 0\,.
\end{aligned}
\end{equation*}
Thus, \eqref{e48f} has been shown. On the other hand, due to \eqref{e58f}, we have
\begin{equation}\label{e59f}
\begin{aligned}
|\zeta(y)-\zeta(x)|&=\zeta(x)\left|e^{\alpha[d(x, x_0)-d(y, x_0)]}-1\right|\\
&  \leq \zeta(x) \alpha | d(x, x_0)- d(y, x_0)| e^{\alpha | d(x, x_0)- d(y, x_0)|}  \\
&\leq \zeta(x) \alpha d(x,y) e^{\alpha d(x,y)}\quad \text{ for all }\; x,y\in G\,.
\end{aligned}
\end{equation}
By \eqref{e59f}, \eqref{e55f}, \eqref{e7f}, \eqref{e1f} and  \eqref{h2b},
\begin{equation*}
\begin{aligned}
\big| I_R^{(2)}\big| &\leq \frac{\alpha}{\delta R} \sum_{x, y\in G} \zeta(x) d^2(x,y) e^{\alpha d(x,y)} |u(x)|^p\omega(x,y)\\
&\leq \frac{ \alpha e^{\alpha s}}{\delta R}  \sum_{x\in G} \zeta(x)|u(x)|^p \sum_{y\in G} \omega(x,y) d^2(x,y)\\
&\leq \frac{  \alpha e^{\alpha s}}{\delta R}  \sum_{x\in G} \zeta(x)|u(x)|^p\mu(x) \leq \frac{ \overline{C} \alpha e^{\alpha s}}{\delta R}  \mathop{\longrightarrow}_{R\to +\infty} 0\,.
\end{aligned}
\end{equation*}
Thus, \eqref{e49f} has been shown. This completes the proof.
\end{proof}

\bigskip
\bigskip

\noindent{\bf Acknowledgement} The authors thank very much the anonymous Referee for his valuable comments, concerning in particular Remarks \ref{oss4}-\ref{oss6}.



\end{document}